\documentclass[11pt]{article}
\usepackage{amsmath, amsfonts, amsthm, latexsym, amssymb}

\def \EE {\mathbb{E}}
\def \RR {\mathbb{R}}
\def \NN {\mathbb{N}}
\def \e  {\epsilon}

\def \bx {{\bf x}}

\def \X  {\mathcal{X}}
\def \S  {\mathrm{sq}}
\def \P  {\mathcal{P}}
\def \N  {\mathcal{N}}
\def \M  {\mathcal{M}}

\def \bB {\mathbf{B}}

\def \oon {\frac{1}{n}}

\def \olp  {\overline{p}}
\def \olen  {\overline{\epsilon_n}}
\def \oleni  {\overline{\epsilon_{ni}}}
\def \olenj  {\overline{\epsilon_{nj}}}

\newtheorem{proposition}{Proposition}
\newtheorem{definition}{Definition}

\title{ Asymptotic Accuracy of the Jackknife Variance Estimator for Certain
Smooth Statistics}
\author{Alex D. Gottlieb}

\date{}

\begin{document}
\maketitle

\begin{abstract}

We show that that the jackknife variance estimator $v_{jack}$ and
the the infinitesimal jackknife variance estimator are
asymptotically equivalent if the functional of interest is a
smooth function of the mean or a smooth trimmed L-statistic.   We
calculate the asymptotic variance of $v_{jack}$ for these
functionals.

\end{abstract}

\section{Introduction}

Let $p$ be a probability measure on a sample space $\X$. Given $n$
samples from $\X$, sampled independently under the probability law
$p$, one desires to estimate the value $T(p)$ of some real
functional $T$ on the space $\P(\X)$ of all probability measures
on $\X$.  Denote by $\e_n$ the map that converts $n$ data points
$x_1,x_2,\ldots,x_n$ into the empirical measure
\begin{equation}
\label{e}
       \e_n(x_1,x_2,\ldots,x_n) \ = \ \oon \sum_{i=1}^n \delta(x_i)
\end{equation}
where $\delta(x_i)$ denotes a point-mass at $x_i$.  The {\it
plug-in estimate} of $T(p)$ given the data $\bx =
(x_1,\ldots,x_n)$ is
\begin{equation}
\label{plug-in}
         T_n \ = \  T(\e_n(\bx)).
\end{equation}
Suppose $T_n$ is an asymptotically normal estimator of $T(p)$, so
that the distribution of $n^{1/2}(T_n - T(p))$ tends to
$\N(0,\sigma^2)$. The jackknife is a computational technique for
estimating $\sigma^2$: one transforms the $n$ original data points
into $n$ pseudovalues and computes the sample variance of those
pseudovalues.

Given the data $\bx = x_1,x_2,\ldots,x_n$, the {\it jackknife
pseudovalues} are
\[
         Q_{ni} \ = \ n T_n(\e_n) \ - \
(n-1)T(\e_{ni}) \qquad \qquad  i = 1,2,\ldots,n
\]
with $\e_n$ as in (\ref{e}) and
\begin{equation}
\label{e-ni}
     \e_{ni} \ = \ \frac{1}{n-1} \sum_{j\ne i} \delta(x_j).
\end{equation}
 The {\it jackknife variance
estimator} is
\begin{equation}
\label{JackVar}
           v_{jack}(x_1,x_2,\ldots,x_n) \ = \ \frac{1}{n-1}\sum_{i=1}^n
                     \left( Q_{ni} - \overline{Q_n} \ \right)^2
\end{equation}
where $\overline{Q_n} = \oon \sum Q_{nj}$.  The variance estimator
$v_{jack}$ is said to be {\it consistent} if
$v_{jack}\longrightarrow \sigma^2$ almost surely as $n \rightarrow
\infty$.  Sufficient conditions for the consistency of $v_{jack}$
are given in terms of the functional differentiability of $T$. An
early result of this kind states that $v_{jack}$ is consistent if
$T$ is strongly Fr\'echet differentiable [Parr85], and it is now
known that $v_{jack}$ is consistent even if $T$ is only
continuously G\^ateaux differentiable as in
Definition~\ref{contGat} below [ST95].

A functional derivative of $T$ at $p$, denoted $\partial T_p$, is
a linear functional that best approximates the behavior of $T$
near $p$ in some sense.  For instance, a functional $T$ on the
space of bounded signed measures $\M(\X)$ is {\it G\^ateaux
differentiable} at $p$ if there exists a continuous linear
functional $\partial T_p$ on $\M(X)$ such that
\[
    \lim_{t \rightarrow 0} \big| t^{-1}\left(T(p + t m) - T(p)\right) \ - \   \partial T_p(m)
    \big| \ = \ 0
\]
for all $m \in \M(\X)$.   More relevant to mathematical statistics
is the concept of Hadamard differentiability, for the fluctuations
of $T(\e_n)$ about $T(p)$ are asymptotically normal if $T$ is
Hadamard differentiable at $p$.  A functional $T:\P(\RR)
\longrightarrow \RR$ is {\it Hadamard differentiable} at $p$ if
there exists a continuous linear functional $\partial T_p$ on
$\M(\RR)$ such that
\[
    \lim_{t\rightarrow 0} \big|t^{-1} \left(T(p + t m_t) - T(p)\right)  \ - \   \partial T_p(m)
    \big| \ = \ 0
\]
whenever $\{m_t\}_{t \in \RR}$ is such that
$\lim\limits_{t\rightarrow 0} m_t = m$ and $m_t(\RR)=0$ for all
$t$, the topology on $\M(\RR)$ being the one induced by the norm
$\|m\| = \sup\limits_{t \in
\RR}\left\{\big|m((-\infty,t])\big|\right\}$.    If $T$ is
Hadamard differentiable at $p$,
 the variance of $n^{1/2}T(\e_n)$ tends to
\begin{equation}
\label{sigmasquared}
           \sigma^2 \ = \ \EE_p\phi_p^2
\end{equation}
as $n \longrightarrow \infty$, where $\phi_p(x)$ is the {\it
influence function}
\begin{equation}
\label{influence}
          \phi_p(x) \ = \  \partial T_p (\delta(x)\ - \ p)
\end{equation}
(this can be shown via the Delta method [vdW98] using Donsker's
theorem).

If $T$ is smooth enough then
$
       n^{1/2}\big( v_{jack} \ - \ \sigma^2 \big)
$ is also asymptotically normal.   In this note we calculate the
asymptotic variance of $v_{jack}$ (i.e., the limit as $n
\longrightarrow \infty$ of the variance of $n^{1/2} v_{jack} $)
for two very well behaved functionals $T$: smooth functions of the
mean $ T(p) \ = \ g \left(\overline{p}\right)$ and smooth trimmed
L-functionals.  In these cases, the asymptotic variance of
$v_{jack}$ equals that of $   \EE_{\e_n} \phi_{\e_n}^2 $, the
estimator of $\sigma^2$ obtained from (\ref{sigmasquared}) by
substituting the empirical measure for $p$.  This is known as the
{\it infinitesimal jackknife} estimator [ST95, p 48]. We are
tempted to conjecture that $v_{jack}$ and the infinitesimal
jackknife variance estimator are asymptotically equivalent for
sufficiently regular functionals $T$, but we have no general
results in this direction.

The literature does not address the accuracy of $v_{jack}$
adequately. In fact, [ST95, Section 2.2.3] gets it wrong,
conjecturing that the asymptotic variance of $ v_{jack} $ should
equal $\hbox{Var}\ \phi_p^2$ for sufficiently regular functionals!
However, Theorem 2 of [Ber84] does contain a general formula for
the variance of $v_{jack}$ which is valid when the functional $T$
has a kind of second-order functional derivative. The theorem
there applies to the trimmed L-functionals we discuss in
Section~\ref{L}, and to many other functionals besides, but it is
hampered by the hypothesis that $p$ have bounded support. We
recommend Theorem~2 of [Ber84] for its generality and its
revelation of the role of second-order differentiability, but our
particular results cannot be derived from it directly.

The text [ST95, p 43] purports to prove that the asymptotic
variance of $n^{1/2} \left( v_{jack} - \sigma^2 \right)$ equals
$\hbox{Var}\ \phi_p^2$ when $T$ is of the form
(\ref{LFunctional}), but there is a mistake there.  We paraphrase
the following definition from [ST95, p 43]:  {\it For probability
measures $p$ and $q$ on the line, let $\rho(p,q)$ denote the
$L^{\infty}$ distance between the cdf's of $p$ and $q$. A
functional $T:\P(\RR) \longrightarrow \RR$ is}
 $\rho$-{\bf Lipschitz differentiable at} $q$
 {\it if
\begin{equation}
\label{LipschitzDiff}
    T(p_k) - T(q_k) \ - \ \partial T_q(p_k-q_k) \ = \
    O\left(\rho(p_k,q_k)^2 \right)
\end{equation}
for all sequences $\{p_k\}$ and $\{q_k\}$ such that $\rho(p_k,q)$
and $\rho(q_k,q)$ converge to $0$}. Assuming that $\hbox{Var}\
\phi_p^2 < \infty$ and $T$ is $\rho$-Lipschitz differentiable, the
authors prove (correctly) that $n^{1/2}\big( v_{jack} - \sigma^2
\big)$ is asymptotically normal with variance $\hbox{Var}\
\phi_p^2$.  They go on to assert that smooth trimmed L-functionals
are $\rho$-Lipschitz differentiable, but this is false (it is not
difficult to construct counterexamples).

A close look at the definition of $\rho$-Lipschitz
differentiability leads one to wonder whether there are any
functionals (besides trivial, linear ones) that satisfy the
definition.  The problem is that $q$ appears on the left hand side
of (\ref{LipschitzDiff}) but not on the right; it is easy to
imagine $p_k$ and $q_k$ that are close to one another in the
$\rho$ metric, yet far enough from $q$ that $\partial
T_q(p_k-q_k)$ badly approximates $T(p_k) - T(q_k)$.  Replacing
$\partial T_q(p_k-q_k)$ by $\partial T_{q_k}(p_k-q_k)$ in the
left-hand-side of (\ref{LipschitzDiff}) might result in a more
useful characteristic of smoothness for a functional $T$. Indeed,
it was this observation that guided our calculations in
Sections~\ref{meanie} and \ref{L}.

In this note we work with modified pseudovalues
\begin{equation}
\label{pseudovalue}
 Q'_{ni}(x_1,x_2,\ldots,x_n) \ = \ (n -1)\left[ T(\e_n) -
T(\e_{ni})\right].
\end{equation}
Substituting $Q'_{ni}$ for $Q_{ni}$ and $\overline{Q_n'} = \oon
\sum Q'_{nj}$ for $\overline{Q_n} = \oon \sum Q_{nj}$ in
(\ref{JackVar}) does not change the value of $v_{jack}$, so one
may compute $v_{jack}$ by the same formula using the $Q'_{ni}$.
Using the modified pseudovalues $Q'_{ni}$ makes it easier to take
advantage of the magic formula $
         (n-1)\left( \e_n - \e_{ni}\right) \ = \ \delta_{x_i} - \e_n
$.

\section{Using pseudovalues to estimate the variance of $\phi_p^2$ }

One aim of this letter is to emphasize that $\hbox{Var}\ \phi_p^2$
is typically {\it not} the asymptotic variance of $n^{1/2}\big(
v_{jack} \ - \ \sigma^2 \big)$, contrary to the assertion of
[ST95, p 42].   However, should one desire an estimate of
$\hbox{Var}\ \phi_p^2$ for some reason, the pseudovalues can be
used to this end.  Once one has already computed $v_{jack}$, the
variance of $\phi_p^2$ is easy to estimate with very little
additional labor: just compute the sample variance of the squares
of the pseudovalues. We prove this, assuming that the functional
$T$ is {\it continuously G\^ateaux differentiable} and $\phi_p$ is
bounded (trimmed L-functionals satisfy these requirements, for
instance).  This section is an interlude whose results will not be
invoked in Sections~\ref{meanie} and \ref{L}, the main part of
this note.

Continuous G\^ateaux differentiability is introduced in [ST95] as
a sufficient condition for the strong consistency of the jackknife
variance estimator.
\begin{definition}
\label{contGat} A functional $T$ is {\bf continuously G\^ateaux
differentiable} at $p$ if it has G\^ateaux derivative $\partial
T_p$ at $p$ and if
\begin{equation}
\label{contGatEq}
 \lim_{k \rightarrow \infty} \sup_{x \in \RR}\left\{ \Big|
           \frac{T(p_k + t_k (\delta(x)-p_k)) - T(p_k)}{t_k} \ - \
               \partial T_p(\delta(x)-p_k) \Big|\right\} \
              \ = \   \ 0
\end{equation}
for any sequence of probability measures $p_k$ whose cdf's
converge uniformly to that of $p$ and for any sequence of real
numbers $t_k$ that converges to $0$.
\end{definition}
The proof in [ST95] that continuous G\^ateaux differentiability
implies strong consistency of the jackknife [ST95, Theorem~2.3]
also serves to prove the following proposition.
\begin{proposition}
\label{P1}
 Suppose that $T:\P(\RR) \longrightarrow \RR$
is continuously G\^ateaux differentiable at $p$, with influence
function $
                 \phi_p(x) \ = \   \partial T_p(\delta(x) - p)
$ satisfying
\[
                 \int |\phi_p(x) | p(dx) < \infty \qquad \int \phi_p(x)p(dx) =
0.
\]
If the data $X_1, X_2, X_3,\ldots$ are iid $p$ then the empirical
measures of the jackknife pseudovalues obtained from the data
converge almost surely to $p\circ \phi_p^{-1}$:
\[
   \e_n(Q'_{n1},Q'_{n2},\ldots,Q'_{nn}) \ \longrightarrow \ p\circ \phi_p^{-1}
\qquad
   \mathrm{a.s.}
\]
\end{proposition}

\noindent {\bf Proof}:  \qquad  Omitted, but cf. the proof of
Theorem~2.3 in [ST95]. \hfill $\square$

Now, suppose that $T:\P(\RR) \longrightarrow \RR$ has a bounded
influence function and satisfies the conditions of
Proposition~\ref{P1}.    Given iid $p$ data $
         X_1,\ X_2,\ldots,\ X_n
$ compute the jackknife pseudovalues
\[
         Q'_{n,1},\  Q'_{n,2},\ldots,\  Q'_{n,n}
\]
and the jackknife estimate $v_{jack}$ based on these pseudovalues.
Set
\[
 \S(x) =
\min\{x^2,\|\phi_p\|^2_{\infty}\} ,
\]
and
\[
                \tau^2 \ = \    \oon \sum_{j=1}^n \left(
\S(Q'_{n,j}) -
 \oon \sum \S(Q'_{n,j})\right)^2.
\]
By Proposition~\ref{P1}, the empirical measure of the jackknife
pseudovalues converges almost surely in $\P(\RR)$ to $p\circ
\phi_p^{-1} $. It follows that $  \tau^2
         \longrightarrow \hbox{Var}\ \phi_p^2 $ almost surely.

One may also estimate $\hbox{Var}\ \phi_p^2$ by applying the
bootstrap to the pseudovalues themselves, just as if the
pseudovalues were actually iid.  To bootstrap, sample $n$ times
with replacement from the empirical measure of the pseudovalues
$Q'_{n,1},\ldots,Q'_{n,n}$, to produce a bootstrap sample
\[   Q^*_{n,1},\ Q^*_{n,2},\ldots,Q^*_{n,n} \]
and compute
\begin{equation}
\label{boot}
              \frac{1}{n^{1/2}}\sum_{i=1}^n \left(\S(Q^*_{n,i})
-\S(Q'_{n,i}) \right).
\end{equation}
Given a triangular array of pseudovalues $Q'_{n,j}$ having the
property that $
      \e_n(Q'_{n,1},\ldots,Q'_{n,n})\ \longrightarrow \ p\circ \phi_p^{-1}
$
as $n \longrightarrow \infty$, one may define $Y_{n,i} =
\S(Q^*_{n,i}) - \oon\sum_j \S(Q'_{n,j})$ and apply the
Lindeberg-Feller Central Limit Theorem to the array
$\{Y_{n,i}\}_{n,i}$ to show that (\ref{boot}) converges in
distribution to $\N(0,\hbox{Var}\ \phi_p^2)$.   But
$\e_n(Q'_{n,1},\ldots,Q'_{n,n})$ almost surely converges to
$p\circ \phi_p^{-1}$ by Proposition~\ref{P1}.  It follows that,
almost surely, (\ref{boot}) converges in distribution to
$\N(0,\hbox{Var}\ \phi_p^2)$.

\section{Functions of the mean}
\label{meanie}

When $q$ is a measure, we denote $\int x q(dx)$ by $\overline{q}$
if the integral is defined.  Let $g \in C^1(\RR)$ and let
\[
   T(m) \ = \ g \left(\overline{m}\right)
\]
be defined for all finite signed measures $m$ with finite first
moment. The functional derivative at $m$ of $T$, evaluated at $q$,
is $
     \partial T_m(q) = g'\left(\overline{m}
     \right)\overline{q}
$; the influence function (\ref{influence}) is $\phi_m(x) =
g'\left(\overline{m}
     \right) \left(x - \overline{m}\right)$.
Suppose that $x_1,x_2,\ldots$ are iid $p$, and $p$ has a finite
second moment. Let $T_n$ denote the plug-in estimator defined in
(\ref{plug-in}). Then the asymptotic variance of
$n^{1/2}\left(T_n-T(p)\right)$ is
\begin{equation}
\label{sigsquared}
         \sigma^2 \ = \ g'(\olp)^2\Big\{ \int x^2 p(dx) - \olp^2 \Big\}.
\end{equation}
Let $v_{jack}$ denote the jackknife variance estimator for
$\sigma^2$.

\begin{proposition}
\label{mean1}
 If $g'$ is (globally) H\"older continuous of order
$h > 1/2$ and $p$ has a finite moment of order $2(1+h)$ then
$n^{1/2} ( v_{jack} - \sigma^2 )$ and $n^{1/2} \left(
\EE_{\e_n}\phi_{\e_n}^2 - \sigma^2 \right)$ have the same limit in
distribution, if any.
\end{proposition}

\noindent {\bf Proof}: \qquad  Set $\Delta_{ni} = \left(Q'_{ni}
-\overline{Q_n'}\right) - \phi_{\e_n}(x_i)$ and note that
\[
    v_{jack}  \ = \ \frac{1}{n-1} \sum_{i=1}^n
    \big(Q_{ni}' - \overline{Q_n'}\big)^2  \ = \  \frac{n}{n-1} \left\{  \EE_{\e_n}\phi_{\e_n}^2 \ + \  \oon
\sum_{i=1}^n \phi_{\e_n}(x_i)\Delta_{ni} \ + \
                            \oon \sum_{i=1}^n \Delta_{ni}^2
                            \right\},
\]
whence
\[
    n^{1/2} \left( v_{jack} - \sigma^2 \right)
     \ = \
     n^{1/2} \left( \EE_{\e_n}\phi_{\e_n}^2 - \sigma^2 \right) \ + \ \frac{n^{1/2}}{n-1}\ \EE_{\e_n}\phi_{\e_n}^2
     \nonumber  \ + \
           \frac{n^{3/2}}{n-1}  \left\{  \oon \sum_{i=1}^n \phi_{\e_n}(x_i)\Delta_{ni} \ + \
                            \oon \sum_{i=1}^n \Delta_{ni}^2 \right\}.
\]
To prove that $n^{1/2} ( v_{jack} - \sigma^2 )$ and $n^{1/2}
\left( \EE_{\e_n}\phi_{\e_n}^2 - \sigma^2 \right)$ have the same
limit in distribution (if any) it suffices to show that
\begin{equation}
\label{tendstozero}
    \frac{n^{1/2}}{n-1}\ \EE_{\e_n}\phi_{\e_n}^2  \ + \
           \frac{n^{3/2}}{n-1}  \left(  \oon \sum_{i=1}^n \phi_{\e_n}(x_i)\Delta_{ni} \ + \
                            \oon \sum_{i=1}^n \Delta_{ni}^2 \right)
\end{equation}
converges almost surely to $0$.

Recall the notation $\e_n$ and $\e_{ni}$ of (\ref{e}) and
(\ref{e-ni}).  The first term in (\ref{tendstozero}) converges
almost surely to $0$ since
\[
      \EE_{\e_n}\phi_{\e_n}^2 \ = \
        \oon \sum_{i=1}^n \phi_{\e_n}^2(x_i) \ = \ \oon \sum_{i=1}^n
g'\left(\olen \right)^2(x_i - \olen)^2
\]
converges almost surely to $\sigma^2$.

 To show that the other
terms tend to zero we need a bound on $\Delta_{ni}$.  Since $g$ is
differentiable, $ g \left( \olenj \right) - g \left(\oleni \right)
\ = \ g'\left(\eta_{ji}\right)\left(\olenj - \oleni \right) $ for
some $\eta_{ji}$ between $\oleni$ and $\olenj$, so that
\[
     Q'_{ni} - \overline{Q_n'} \ = \
          \frac{n-1}{n} \sum_{j=1}^n
          \big(  g \left( \olenj \right) - g \left(\oleni \right)\big)
          \ = \
          \frac{n-1}{n} \sum_{j=1}^n
          g'\left(\eta_{ji}\right)\left(\olenj - \oleni \right).
\]
Therefore, since $\phi_{\e_n}(x_i) = g'\left(\olen \right)(x_i -
\olen) = \oon\sum_j g'\left(\olen \right)(x_i - x_j)$,
\begin{eqnarray*}
     \Delta_{ni}
     \ = \
     \left(Q'_{ni} - \overline{Q_n'}\right) - \phi_{\e_n}(x_i)
     & = &
        \frac{n-1}{n} \sum_{j=1}^n
          g'\left(\eta_{ji}\right)\left(\olenj - \oleni \right)
          \  - \  \oon\sum_{j=1}^n g'\left(\olen \right)(x_i - x_j) \\
     &  = &
                   \oon \sum_{j=1}^n
         \left( g'\left(\eta_{ji}\right) - g'\left(\olen \right)\right)(x_i - x_j
         ).
\end{eqnarray*}
But $g'$ is H\"older continuous of order $h$ and
$|\eta_{ji}-\olen|< \max\{|\olenj - \olen|,|\oleni - \olen|\}$, so
\[
            \left|  g'\left(\eta_{ji}\right) - g'\left(\olen
            \right)\right| \ \le \  C \big( |\olenj - \olen|^h +
            |\oleni - \olen|^h \big)
            \ \le \  C (n-1)^{-h} \big( |\olen - x_j|^h +
            |\olen - x_i|^h \big),
\]
where $C$ is a global H\"older constant for $g'$.  It follows that
\[
     \left| \Delta_{ni} \right|
     \ = \
     C (n-1)^{-h} \oon \sum_{j=1}^n   \big( |\olen - x_j|^h +
            |\olen - x_i|^h  \big) \big( |\olen - x_j| +
            |\olen - x_i|  \big).
\]
With this bound on $\Delta_{ni}$, and  assuming that $p$ has a
finite moment of order $2(1+h)$, it may be shown that
\[
     \oon \sum_{i=1}^n \Delta_{ni}^2 \ = \ O_s\big(n^{-2h}),
\]
and then, by the Cauchy-Schwartz inequality, that
\[
         \Big\vert \oon \sum_{i=1}^n \phi_{\e_n}(x_i)\Delta_{ni}
         \Big\vert \ = \ O_s\big( n^{-h} \big).
\]
  The preceding estimates and
the assumption that $h > 1/2$ imply that the last two terms in
(\ref{tendstozero}) converge to almost surely to $0$.  Thus,
$n^{1/2} \left( v_{jack} - \sigma^2 \right)$ and $n^{1/2} \left(
\EE_{\e_n}\phi_{\e_n}^2 - \sigma^2 \right)$ have the same limit in
distribution, if any. \hfill $\square$

If we strengthen the smoothness assumption on $g$ and the moment
assumption on $p$ then we can calculate the limit in distribution
of $n^{1/2}\left( \EE_{\e_n}\phi_{\e_n}^2 - \sigma^2 \right)$.
 Suppose that $g''$ is bounded (so that $g'$ is globally
Lipschitz) and H\"older continuous of order $r > 0$, and suppose
that $p$ has a finite fourth moment.  Then
\[
      \phi_{\e_n}(x_i)
         \ = \
           g'\left(\olen \right) \left( x_i - \olen \right)
         \ = \
           \left[ g'\left( \olp \right) + g''\left(\olp\right)\left( \olen - \olp \right)
        + O_s \big( n^{-(r+1)/2} \big) \right] \left( x_i - \olen
        \right),
\]
so that
\begin{eqnarray*}
   \EE_{\e_n}\phi_{\e_n}^2
    & = &
    \oon \sum_{i=1}^n \phi_{\e_n}^2(x_i)
   \ = \
    \left[ g'\left( \olp \right) + g''\left(\olp\right) \left( \olen - \olp
    \right)\right]^2
    \oon \sum_{i=1}^n  \left( x_i - \olen \right)^2
        \ + \ O_s \big( n^{-(r+1)/2} \big)
    \\
    & = &
     \left[
      g'\left( \olp \right)^2 + 2 g'\left( \olp \right)
      g''\left(\olp\right)\left( \olen - \olp \right)
     \right] \oon \sum_{i=1}^n\left( x_i - \olen \right)^2
       \ + \  O_s \big( n^{-(r+1)/2} \big).
     \\
\end{eqnarray*}
From formula (\ref{sigsquared}) for $\sigma^2$ we see that
\begin{eqnarray}
      n^{1/2}\left( \EE_{\e_n}\phi_{\e_n}^2 - \sigma^2 \right)
      & = &
      g'\left( \olp \right)^2 n^{1/2}\left( \oon \sum_{i=1}^n  \left( x_i - \olen
      \right)^2 \ - \ \Big\{ \int x^2 p(dx) - \olp^2 \Big\}\right)
      \nonumber \\
      &   & \ + \
       2 g'\left( \olp \right) g''\left(\olp\right)
       n^{1/2}\left( \olen - \olp \right) \oon \sum_{i=1}^n\left( x_i - \olen
      \right)^2
      \ + \ O_s \left( n^{-r/2}\right) . \label{fluctuations}
\end{eqnarray}
Set $Z_n = n^{1/2}\left( \olen - \olp \right)$ and
\[
         Y_n \ = \ n^{1/2}\left( \oon \sum_{i=1}^n  \left( x_i - \olen
      \right)^2 \ - \ \Big\{ \int x^2 p(dx) - \olp^2
      \Big\}\right).
\]
Since $p$ has a finite fourth moment, the random vector
$(Y_n,Z_n)$ has a Gaussian limit by the Central Limit Theorem.
Equation (\ref{fluctuations}) shows that $n^{1/2}\left(
\EE_{\e_n}\phi_{\e_n}^2 - \sigma^2 \right)$ is asymptotically
normal with variance $ (a,b) \Gamma (a,b)^{tr} $, where $(a,b) =
\left( g'( \olp )^2,\  2 g'( \olp ) g''(\olp) \right)$ and
$\Gamma$ denotes the asymptotic covariance matrix for $(Y_n,Z_n)$.

In view of Proposition~\ref{mean1}, we find that if $g''$ is
bounded and H\"older continuous of order $r > 0$, and if $p$ has a
finite fourth moment, then the asymptotic variance of $n^{1/2}
\left( v_{jack} - \sigma^2 \right)$ equals $ (a,b) \Gamma
(a,b)^{tr} $.    In contrast, under the same conditions on $p$ and
$g$ it may be shown that $\hbox{Var}\ \phi_p^2 = a^2
\Gamma_{1,1}$.

\section{Trimmed L-statistics}
\label{L}

Suppose that $\ell:(0,1) \longrightarrow \RR$ is supported on
$[\alpha,1-\alpha]$ for some $0 < \alpha < 1/2$, and let
\begin{equation}
\label{LFunctional}
     L(p) \ = \ \int_0^1 P^{-1}(s)\ell(s)ds.
\end{equation}
Here $P^{-1}$ denotes the quantile function for $p$, i.e.,
$
     P^{-1}(s) = \min\{x: P(x) \ge s\}
$ for $0<s<1$ where $P$ denotes the cdf of $p$.  A plug-in
estimate for $L$ is called a {\it trimmed L-statistic}, or a
trimmed {\it linear combination of quantiles}.  (It is called {\it
trimmed} because the restricted support of $\ell$ discards
outliers.) L-statistics are good for robust estimation of a
location parameter.

Now assume that $\ell$ is continuous.  Then $L$ is Hadamard
differentiable (and the L-statistics are asymptotically normal) at
all $p \in \P(\RR)$ [vdW98, Lemma 22.10].    The functional
derivative at $p$ of $L$, evaluated at a bounded signed measure
$m$, is
\[
          \partial L_p(m) \ = \ -\int \ell(P(x))M(x)dx
\]
where $M(x) = m((-\infty,x])$.    The asymptotic variance of the
L-statistics is
\[
       \sigma^2  \ = \ \int\int \ell(P(y))\Gamma(y,z)\ell(P(z))dydz,
\]
where
\begin{equation}
\label{Gamma}
    \Gamma(y,z) \ = \ P(y)\wedge P(z) \ - \ P(y)P(z).
\end{equation}
This formula is obtained via Donsker's Theorem: {\it Let $P_n$
denote the cdf of $\e_n$, a random bounded function. Then $
n^{1/2}(P_n(t)-P(t))$ converges in law to a Gaussian process
$\{\bB(t)\}_{t \in \RR}$ with covariance}
\begin{equation}
         \Gamma(s,t) \ = \  \EE_p\left[\bB(s)\bB(t) \right] \ = \ P(s)\wedge P(t) -
           P(s)P(t).
\label{BB}
\end{equation}
Finally, the influence function is
\begin{equation}
\label{LFunInfluence}
 \phi_p(x) \ = \ \partial L_p(\delta(x) - p) \ = \
-\int \ell(P(y))(H_x - P)(y)dy,
\end{equation}
where $H_{x}$ denotes the cdf of $\delta(x)$.  Note that $
\sigma^2 = \EE_p \phi_p^2 $ and
\[
         \EE_{\e_n}\phi_{\e_n}^2 \ = \ \int\int
         \ell(P_n(y)) \left[  P_n(y)\wedge P_n(z) - P_n(y)P_n(z)
         \right]
\ell(P_n(z))dydz.
\]

Let $v_{jack}$ denote the jackknife variance estimator for
$\sigma^2$.  We find that the $v_{jack}$ is asymptotically
equivalent to $\EE_{\e_n}\phi_{\e_n}^2$ and asymptotically normal:

\begin{proposition}
\label{LFunProp}
 Suppose $p$ has no point masses and $\ell'$ is
H\"older continuous of order $h > 1/2$.  Then
\begin{equation}
\label{equivalent}
     n^{1/2} \left( v_{jack} - \sigma^2 \right) \ = \
     n^{1/2}\left( \EE_{\e_n}\phi_{\e_n}^2 - \sigma^2 \right) \ + \
              O_s\big(n^{1/2-h}\big)
\end{equation}
and converges in law to the Gaussian random variable $ Y + Z$,
where
\begin{eqnarray}
       Y & = &
         \int \int \ell(P(y))
            \left\{ \bB(y\wedge z) - P(y)\bB(z) - \bB(y)P(z) \right\} \ell(P(z))
            dydz  \nonumber \\
       Z & = &
         2  \int \int
               \ell'(P(y))\bB(y)\Gamma(y,z)
               \ell(P(z))dydz  \label{YnZ}
\end{eqnarray}
and $\bB$ denotes the Brownian Bridge (\ref{BB}).
\end{proposition}

\noindent {\bf Proof}:  \qquad  We prove first that $n^{1/2}\left(
\EE_{\e_n}\phi_{\e_n}^2 - \sigma^2 \right)$ converges in law to
$Y+Z$, and afterwards we establish (\ref{equivalent}).

\noindent Define
\begin{eqnarray}
          Y_n & = & n^{1/2}\Big( \sum_{i=1}^n \phi_p(x_i)^2 -  \sigma^2 \Big) \nonumber \\
          Z_n & = & -2n^{-1/2} \sum_{i=1}^n \phi_p(x_i)
           \int \ell'(P(y))(P_n-P)(y)\left( H_{x_i} - P_n \right)(y)
           dy.       \label{YnandZn}
\end{eqnarray}
We claim that $Y_n$ converges in law to $Y$ and $Z_n$ converges in
law to $Z$.  To see this, substitute (\ref{LFunInfluence}) for
$\phi_p$ in the definitions of $Y_n$ and $Z_n$, and apply
Donkser's Theorem.  Substituting (\ref{LFunInfluence}) for
$\phi_p$ yields
\begin{eqnarray*}
    Y_n
    & = &   \int \int \ell(P(y))
            n^{1/2}\Big(\oon \sum_{i=1}^n H_{x_i}(y)H_{x_i}(z) \ - \
P(y)\wedge P(z)\Big) \ell(P(z)) dydz \\
    &   & \ - \ \int \int \ell(P(y)) P(y)n^{1/2}(P_n-P)(z) \ell(P(z)) dydz \\
    &   & \qquad  \ - \ \int \int \ell(P(y)) n^{1/2}(P_n-P)(y)P(z) \ell(P(z))
    dydz \\
     Z_n & = &  2 n^{-1/2}
   \sum_{i=1}^n \int \int \ell'(P(y))(P_n-P)(y)\left( H_{x_i} - P_n \right)(y)
                   \ell(P(z))(H_{x_i} - P)(z) dydz  \\
         & = & 2  \int \int
               \ell'(P(y))n^{1/2}(P_n-P)(y)\Big( \oon \sum_{i=1}^n
H_{x_i}(y)H_{x_i}(z) - P_n
               (y)P_n(z) \Big) \ell(P(z))dydz.
\end{eqnarray*}
Note that $ \oon \sum H_{x_i}(y)H_{x_i}(z) - P_n(y)P_n(z)$ in the
expression for $Z_n$ converges almost surely to $\Gamma(y,z)$ of
(\ref{Gamma}).   Also, in the expression for $Y_n$,
\[
      n^{1/2}\Big(\oon \sum_{i=1}^n H_{x_i}(y)H_{x_i}(z) -
P(y)\wedge P(z)\Big) \ = \ n^{1/2}\left( P_n(y)\wedge P_n(z) -
P(y)\wedge P(z) \right)
\]
converges in law to the Gaussian process $ \bB(y\wedge z)$.
 Writing $M_{ni} =
H_{x_i} - P_n$, we find that
\begin{eqnarray}
      \phi_{\e_n}(x_i)
         & = &  - \int \left\{ \ell(P(y) + \ell'(P(y))(P_n - P)(y) +
                   O_s(n^{-h})\right\} M_{ni}(y) dy  \nonumber \\
       & = &    \phi_p(x_i) \ - \  \int \ell'(P(y))(P_n - P)(y)M_{ni}(y)
       dy     \  + \ O_s\big(n^{-h}\big). \label{parTy}
\end{eqnarray}
Equations  (\ref{parTy}) and (\ref{YnandZn}) imply that
\begin{eqnarray*}
    n^{1/2} \left(\EE_{\e_n}\phi_{\e_n}^2 -
    \sigma^2 \right) \ = \   Y_n \ + \ Z_n
    & + &
     n^{-1/2}
    \sum_{i=1}^n \left( \int \ell'(P(y))(P_n(y) - P(y))M_{ni}(y) dy \right)^2
      \\
   & + &  O_s\big(n^{1/2-h}\big).
\end{eqnarray*}
But the third term on the right hand side of the last equation is
$O_s\big(n^{-1/2}\big)$, since
\begin{eqnarray*}
    &&  \oon \sum_{i=1}^n \left( \int \ell'(P(y))(P_n(y) - P(y))M_{ni}(y) dy
\right)^2 \\
     && \ = \ \oon \sum_{i=1}^n \int\int \ell'(P(y))(P_n - P)(y)
     \ell'(P(z))(P_n - P)(z)M_{ni}(y) M_{ni}(z) dydz \\
     && \ = \ \int\int \ell'(P(y))(P_n - P)(y)
     \ell'(P(z))(P_n - P)(z)\oon \sum_{i=1}^n M_{ni}(y) M_{ni}(z) dydz \\
\end{eqnarray*}
and
\[
       \oon \sum_{i=1}^n M_{ni}(y) M_{ni}(z) \ = \ \oon
       \sum_{i=1}^n H_{x_i}(y)H_{x_i}(z) \ - \ P_n(y)P_n(z),
\]
converges almost surely to $\Gamma(y,z)$.  Thus,
\[
 n^{1/2} \left( \EE_{\e_n}\phi_{\e_n}^2 - \sigma^2 \right) \ = \ Y_n \ + \ Z_n
 \ + \ O_s\big(n^{-h}\big),
\]
so that $n^{1/2} \left( \EE_{\e_n}\phi_{\e_n}^2 - \sigma^2
\right)$ converges in law to   $ Y + Z$, a Gaussian random
variable.

It remains to establish (\ref{equivalent}).  To this end it
suffices to show that
\begin{equation}
\label{suffices}
       \max_{1\le i \le n} \left\{ \left| Q_{ni}' -
       \overline{Q_n'} - \phi_{\e_{n}}(x_i) \right| \right\} = O_s\big(
       n^{-h}\big),
\end{equation}
for then, since $v_{jack} = (n-1)^{-1} \sum
    (Q_{ni}' - \overline{Q_n'})^2$, it would follow that
\begin{eqnarray*}
     n^{1/2} \left( v_{jack} -
    \sigma^2 \right) & = &
     n^{1/2}\Big( \frac{1}{n-1} \sum_{i=1}^n \phi_{\e_n}^2(x_i)
           \ - \  \sigma^2\Big)  \ + \ O_s\big(n^{1/2-h}\big) \nonumber \\
    & = &   n^{1/2}\left( \EE_{\e_n}\phi_{\e_n}^2 - \sigma^2 \right) \ + \
              O_s\big(n^{1/2-h}\big).
\end{eqnarray*}
 Let $P_{ni}$ denote the cdf of
$\e_{ni}$.  Integration by parts of (\ref{LFunInfluence}) shows
that
\begin{equation}
\label{Influence2}
    \phi_{\e_n}(x_i) \ = \ \int x d \left[ \ell(P_n)(
    H_{x_i} - P_n)(y)\right]
\end{equation}
(the boundary term vanishes because (\ref{LFunctional}) is
trimmed).  Suppose $x_1,x_2,x_3,\ldots$ are distinct (we are
assuming that $p$ has no point masses, so this is the case
 almost surely).  Then (\ref{Influence2}) becomes
\begin{eqnarray*}
    \phi_{\e_n}(x_i)
    & = &
     x_i \ell (P_n(x_i)) \ +  \  \sum_{j:\ x_j > x_i} x_j \left\{ \ell(P_n(x_j))
    \ - \  \ell\big(P_n(x_j)- 1/n \big) \right\} \\
    &   &
     - \ \sum_{j=1}^n x_j \left\{ \ell(P_n(x_j)) P_n(x_j)
    \ - \  \ell\big(P_n(x_j)- 1/n \big) \big(P_n(x_j)- 1/n\big) \right\},
\end{eqnarray*}
which we rewrite as $\phi_{\e_n}(x_i) = A + B_i + C_i + D_i \ $
with
\begin{eqnarray}
    A & = &
     - \ \oon \sum_{j=1}^n x_j \ell\big(P_n(x_j)- 1/n \big)  \nonumber \\
    B_i & = &
     x_i \ell (P_n(x_i)) \nonumber \\
    C_i & = &
    - \ \sum_{j:x_j \le x_i} x_j \left\{ \ell(P_n(x_j)) -  \ell(P_n(x_j)- 1/n )\right\} P_n(x_j)   \nonumber \\
    D_i & = &
     \sum_{j:x_j > x_i} x_j \left\{ \ell(P_n(x_j)) -  \ell(P_n(x_j)- 1/n )\right\} \left( 1 -P_n(x_j)\right) .
    \label{ABCD}
\end{eqnarray}
For $1 \le i \le n$, let
\[
         \zeta_{ni}(x) \ = \ (n-1) \int_{P_{ni}(x) - \frac{1}{n-1}}^{P_{ni}(x)}  \ell(s)
         ds.
\]
Observe that $\ell\big(\zeta_{ni}(x)\big) =
\ell\big(\zeta_{nk}(x)\big)$ if $x < \min\{x_i, x_k\}$ or if $x
> \max\{x_i, x_k\}$, and
\begin{eqnarray}
  \zeta_{nk}(x) - \zeta_{ni}(x) \ = \
       (n-1)\int_{P_{ni}(x)}^{P_{ni}(x) + \frac{1}{n-1}}
       \ell(s)-\ell\big(s - 1/(n-1)\big)ds
       &\quad \hbox{if} &
        x_i < x < x_k   \nonumber \\
  \zeta_{nk}(x) - \zeta_{ni}(x) \ = \
       -(n-1)\int_{P_{ni}(x) - \frac{1}{n-1}}^{P_{ni}(x)}
       \ell(s)-\ell\big(s - 1/(n-1)\big)ds
       & \quad \hbox{if} &
        x_k < x < x_i.
       \label{observation}
\end{eqnarray}
Thus $L(\e_{ni}) = \frac{1}{n-1}\sum\limits_{j:j\ne i} x_j
 \zeta_{ni}(x_j)$ and
\begin{eqnarray*}
   Q'_{ni} - \overline{Q_n'} & = &
           -\sum_{j:j \ne i} x_j \zeta_{ni}(x_j)
           \ + \ \oon \sum_{k=1}^n \sum_{j:j \ne k}
                                          x_j \zeta_{nk}(x_j)  \\
       & = &
             -\ \oon \sum_{k=1}^n x_k \zeta_{ni}(x_k)
               \ + \
               \oon \sum_{k=1}^n x_i \zeta_{nk}(x_i)
               \  + \ \oon \sum_{k=1}^n
       \sum_{j:j \ne k,i} x_j \left\{ \zeta_{nk}(x_j)-
       \zeta_{ni}(x_j)\right\}    \\
     & = &
         -\ \oon \sum_{k=1}^n x_k \zeta_{ni}(x_k)
               \ + \
               \oon \sum_{k=1}^n x_i \zeta_{nk}(x_i)
        \  + \  \oon \sum_{j:x_j < x_i} \  \sum_{k:x_k < x_j}
         x_j \left( \zeta_{nk}(x_j)-
       \zeta_{ni}(x_j) \right) \\
      &   &
        \  + \     \oon \sum_{j:x_j > x_i} \ \sum_{k:x_k > x_j}
         x_j \left( \zeta_{nk}(x_j)-
       \zeta_{ni}(x_j) \right).
\end{eqnarray*}
Using (\ref{observation}) we find that $Q'_{ni} - \overline{Q_n'}
= A_i'+B_i'+C_i'+D_i'\ $ with
\begin{eqnarray}
    A_i' & = &
                -\ \oon \sum_{j=1}^n x_j \zeta_{ni}(x_j) \nonumber \\
    B_i' & = &
               \oon \sum_{j=1}^n x_i \zeta_{nj}(x_i)
               \nonumber \\
    C_i' & = &
   (n-1) \sum_{j:x_j < x_i}
         x_j \left(P_n(x_j) - 1/n\right) \int_{P_{ni}(x_j) - \frac{1}{n-1}}^{P_{ni}(x_j)}
       \ell(s)-\ell\big(s - 1/(n-1)\big)ds  \nonumber \\
    D_i' & = &
      (n-1) \sum_{j:x_j > x_i} x_j \left(1 - P_n(x_j)\right)
         \int_{P_{ni}(x_j)}^{P_{ni}(x_j) + \frac{1}{n-1}}
       \ell(s)-\ell\big(s - 1/(n-1)\big)ds.
    \label{A'B'C'D'}
\end{eqnarray}
The sequence $\{P_n\}$ converges almost surely to $P$ and hence it
is almost surely tight.  Thus there exists a (random) bound $M>0$
such that $P_n(x) < \alpha/2$ if $x<M$ and $P_n(x) > 1 - \alpha/2$
if $x > M$.   Since $\ell$ vanishes off of $[\alpha,1-\alpha]$, it
follows that $B_i = 0$ if $|x_i| > M$, and $B'_i = 0$ if $|x_i| >
M$ and $1/(n-1)<\alpha/4$.   Similarly, if $n$ is sufficiently
large, the sums defining $A',C',D',A_i',C_i'$ and $D_i'$ in
(\ref{ABCD}) and (\ref{A'B'C'D'}) may be replaced with sums over
$j$ such that $|x_i| > M$.  Thus
\begin{eqnarray*}
     |A'_i - A|
     & \le &
     M \frac{n-1}{n} \sum_{j=1}^n  \int_{P_{ni}(x_j) - \frac{1}{n-1}}^{P_{ni}(x_j)}
                \left|  \ell(s) - \ell\big(P_n(x_j)- 1/n \big) \right| ds
      \\
     |B'_i - B_i|
     & \le &
     M \frac{n-1}{n} \sum_{j=1}^n
      \int_{P_{nj}(x_i) - \frac{1}{n-1}}^{P_{nj}(x_i)} \left| \ell(s) - \ell (P_n(x_i))
     \right|ds
\end{eqnarray*}
are both $O_s(1/n)$ since $\ell$ is differentiable.   For $n > 1
\NN$ and $s \in [1/n,1]$, let $t_n(s)$ be a number between $s
-1/n$ and $s$ such that $
        \ell'(t_n(s)) = n \left( \ell(s)-\ell\big(s -
         1/n \big)\right) $.
 (The functions $t_n$ may be chosen to be continuous, since $\ell'$ is
 continuous.)
We now have
\begin{eqnarray*}
     |C'_i - C_i|
     & \le &
            M \big| \ell'(t_n(P_n(x_i)))\big| P_n(x_i)
      \   +  \
      \frac{M}{n} \sum_{j:x_j < x_i}
       \int_{P_{ni}(x_j) - \frac{1}{n-1}}^{P_{ni}(x_j)}
       \big| \ell'(t_{n-1}(s)) \big| ds
       \\
       &    &
       + \
        M \sum_{j:x_j < x_i}
         P_n(x_j) \int_{P_{ni}(x_j) - \frac{1}{n-1}}^{P_{ni}(x_j)}
       \Big| \ell'(t_{n-1}(s)) - \ell'(t_n(P_n(x_j))) \Big|
       ds
       \\
       |D'_i - D_i|
       &  \le &
       M \sum_{j:x_j > x_i}  \left(1 - P_n(x_j)\right)
         \int_{P_{ni}(x_j)}^{P_{ni}(x_j) + \frac{1}{n-1}}
       \Big|  \ell'(t_{n-1}(s)) - \ell'(t_n(P_n(x_j))) \Big|ds.
\end{eqnarray*}
But $ \ell'(t_{n-1}(s)) - \ell'(t_n(P_n(x_j))) = O\big(
n^{-h}\big)$ throughout the interval of integration because of the
H\"older continuity of $\ell'$, and so $|C'_i - C_i|$ and $|D'_i -
D_i|$ are both $O_s\big( n^{-h}\big)$ uniformly in $i$. The
preceding estimates show that
\[
        \left| Q_{ni}' -
       \overline{Q_n'} - \phi_{\e_{n}}(x_i) \right|
       \ \le \
       |A'_i - A| + |B'_i - B| + |C'_i - C_i| + |D'_i - D_i|
       \ = \
       O_s\big( n^{-h}\big)
\]
uniformly in $i$, establishing (\ref{suffices}).     \hfill
$\square$

Proposition~\ref{LFunProp} is also true as stated for $ L(p)= \int
x \ell(P(x))p(dx) $, which is not exactly the same as the
L-functional (\ref{LFunctional}) but has the same functional
derivative.   An argument similar to the one above shows that the
asymptotic variance of $n^{1/2} \left( v_{jack} - \sigma^2
\right)$ equals $\hbox{Var}\ (Y+Z)$ with $Y$ and $Z$ as in
(\ref{YnZ}).   On the other hand, one can show that $\hbox{Var}\
\phi_p^2 = \hbox{Var}\ Y$.  This is contrary to [ST95, p 43],
where it is asserted that $\hbox{Var}\ Y$ is the asymptotic
variance of $n^{1/2} \left( v_{jack} - \sigma^2 \right)$.

\section{Acknowledgments}

Thanks to Steve Evans for his advice and encouragement. Thanks to
Rudolf Beran.  The author is supported by the Austrian START
project {\it Nonlinear Schr\"odinger and quantum Boltzmann
equations}.

\section{References}

\noindent [Ber84] \quad R. Beran.  Jackknife approximations to
bootstrap estimates.

\qquad  \qquad {\it The Annals of Statistics} 12 (1): 101 - 118,
1984.

\noindent [Parr85] \quad W. Parr.  Jackknifing differentiable
statistical functions.

\qquad  \qquad {\it Journal of the Royal Statistical Society B} 47
(1): 56 - 66, 1985.

\noindent [ST95] \quad   J. Shao and D. Tu.  {\it The Jackknife
and Bootstrap}. Springer-Verlag, New York, 1995.

\noindent [vdW98] \quad   A.W. van der Waart.  {\it Asymptotic
Statistics}. Cambridge University Press, 1998.

\end{document}